\documentclass[12pt,a4paper,twoside]{article}

\newcommand{\versiondate}{30 April 2004}

\usepackage{amsmath,amssymb,amscd,amsthm,latexsym} 
\usepackage[all]{xy}

\theoremstyle{plain}  

\newtheorem{theorem}{Theorem}[section]
\newtheorem*{theorem*}{Theorem}

\newtheorem{corollary}[theorem]{Corollary}
\newtheorem{lemma}[theorem]{Lemma}
\newtheorem{proposition}[theorem]{Proposition}

\theoremstyle{definition}

\newtheorem{definition}[theorem]{Definition}

\theoremstyle{remark}

\newtheorem{example}[theorem]{Example}
\newtheorem{remark}[theorem]{Remark}

\newtheorem*{claim*}{Claim}

\numberwithin{equation}{section}

\renewcommand{\geq}{\geqslant}

\newcommand{\lto}{\longrightarrow}
\newcommand{\lmapsto}{\longmapsto}

\newcommand{\HH}{\mathbb{H}}

\newcommand{\ftil}{\widetilde{f}}
\newcommand{\into}{\hookrightarrow}

\newcommand{\kk}{\mathbf{k}}
\newcommand{\blank}{\text{--}}
\newcommand{\OO}{\mathcal{O}}
\newcommand{\Hom}{\mathrm{Hom}}

\newcommand{\Ext}{\mathrm{Ext}}
\newcommand{\shom}{\mathcal{H}\hspace{-1pt}\mathit{om}}

\DeclareMathOperator{\im}{im}

\begin{document}

\null
\vskip 2em
\begin{center}
  {\LARGE Homological algebra of twisted quiver bundles \par}
\vskip 1.5em
  {\large Peter B. Gothen\footnotemark[1]$^{,}$\footnotemark[2] 
   and Alastair D. King\footnotemark[1] \par
\vskip 1em
   \versiondate}
\end{center}  
  \par
  \vskip 1.5em  

\footnotetext[1]{
Members of VBAC (Vector Bundles on
Algebraic Curves), which is partially supported by EAGER (EC FP5
Contract no.\ HPRN-CT-2000-00099) and by EDGE (EC FP5 Contract no.\ 
HPRN-CT-2000-00101).}

\footnotetext[2]{
Partially supported by the
Funda{\c c}{\~a}o para a Ci{\^e}ncia e a Tecnologia (Portugal) through
the Centro de Matem{\'a}tica da Universidade do Porto and through
grant no.\ SFRH/BPD/1606/2000.}

\begin{abstract}
  Several important cases of vector bundles with extra structure 
  (such as Higgs bundles and triples) may be regarded as examples
  of twisted representations of a finite quiver in the category
  of sheaves of modules on a variety/manifold/ringed space.
  We show that the category of such representations
  is an abelian category with enough injectives
  by constructing an explicit injective resolution.
  Using this explicit resolution, we find a long exact
  sequence that computes the $\Ext$ groups in this new category
  in terms of the $\Ext$ groups in the old category.
  The quiver formulation is directly reflected in the form of the 
  long exact sequence.
  We also show that under suitable circumstances, 
  the $\Ext$ groups are isomorphic to certain hypercohomology groups.
\end{abstract}

\section{Introduction}
\label{sec:introduction}

A \emph{quiver} $Q$ is a directed graph specified by a set of vertices
$Q_0$, a set of arrows $Q_1$ and head and tail maps
\begin{displaymath}
  h,t \colon Q_1 \lto Q_0.
\end{displaymath}
Throughout the paper we shall assume that $Q$ is finite, i.e.\ the
sets $Q_0$ and $Q_1$ are finite.

Let $\mathcal{A}$ be a category.  A \emph{representation} of
a quiver $Q$ in $\mathcal{A}$ is a cor\-res\-pon\-dence which associates an
object $V_i$ of $\mathcal{A}$ to each vertex $i \in Q_0$ and a
morphism
\begin{displaymath}
  \phi_{a} \colon V_{ta} \to V_{ha}
\end{displaymath}
to each arrow $a \in Q_1$.
There is an obvious notion of a morphism $f \colon (V,\phi) \to 
(W,\psi)$ between
representations of $Q$: it consists of a family of morphisms $f_i
\colon V_i \to W_i$ such that the diagrams
\begin{displaymath}
  \xymatrix{
    V_{ta} \ar[d]^{f_{ta}} \ar[r]^{\phi_{a}} & V_{ha} \ar[d]^{f_{ha}} \\
    W_{ta} \ar[r]^{\psi_{a}}                 & W_{ha} \\
  }
\end{displaymath}
commute for every $a \in Q_1$.  In this way the representations of $Q$
form a category; note that this is an abelian category whenever
$\mathcal{A}$ is.

A fundamental case is when $\mathcal{A}$ is the category of vector
spaces over a field $\kk$.  We shall refer to a representation of $Q$
in this category as a $Q$-vector space or simply as a representation
of $Q$.  There is a related algebraic notion: given a quiver $Q$ one
can construct its so-called path algebra (see
Section~\ref{sec:path-algebra}) and the category of left modules over
the path algebra is equivalent to the category of representations of
$Q$.  This has been a very fruitful idea in the theory of
representations of algebras (see, e.g., Benson \cite{benson:1998}).

Another case of great interest is when $\mathcal{A}$ is the category
of algebraic vector bundles over an an algebraic variety $X$, or the
category of holomorphic vector bundles over a complex manifold $X$.  A
representation of a quiver $Q$ in this category is called a
\emph{quiver bundle} or, for short, a \emph{$Q$-bundle}.  This notion
unifies many of the vector bundles with extra structure which have
been studied previously; we mention some examples.  If $Q$ is the
quiver
\begin{displaymath}
        \xymatrix{
          *!/u .9ex/{\underset{0}{\bullet}}
          & *!/u .9ex/{\underset{1}{\bullet}} \ar[l]_{a}
    }\ ,
\end{displaymath}
then a $Q$-bundle is a holomorphic triple as studied by
Garc{\'\i}a-Prada \cite{garcia-prada:1994} and Bradlow and
Garc{\'\i}a-Prada \cite{bradlow-garcia-prada:1996}.  In the special
case $V_1 = \OO_X$ this becomes equivalent to the notion of Bradlow pair \cite{bradlow:1991}.  More generally, consider the quiver 
\begin{displaymath}
        \xymatrix{
          *!/u .9ex/{\underset{0}{\bullet}}
          & *!/u .9ex/{\underset{1}{\bullet}} \ar[l]_-{a_1}
          & {\cdots} \ar[l]_-{a_2}
          & *!/u .9ex/{\underset{m}{\bullet}} \ar[l]_-{a_m}
    }\ . 
\end{displaymath}
In this case a $Q$-bundle is called a holomorphic chain, these where
studied by \'Alvarez-C\'onsul and Garc{\'\i}a-Prada
\cite{alvarez-garcia-prada:2001}.  There is a generalization of the
notion of $Q$-bundle to that of \emph{twisted $Q$-bundle}, where
the morphisms $\phi_{a}$ are twisted by a vector bundle $M_a$ (see
Section~\ref{sec:Q-sheaves} for precise definitions).
An important example of this is that of Higgs bundles on a curve $X$, 
introduced by Hitchin \cite{hitchin:1987}: 
these are pairs $(E,\phi)$, where $E$ is
a vector bundle and $\phi$ is a morphism $\phi \colon E \lto E \otimes
K_X$, where $K_X$ is the canonical bundle of $X$.  Thus a Higgs bundle
is a twisted $Q$-bundle where $Q$ is given by
\begin{equation}
  \label{eq:higgs-quiver}
    \xymatrix{
      *!/u .9ex/{\underset{0}{\bullet}} \ar@(ur,dr)[]^{a}     
    }\ .  
\end{equation}

In various applications it has been important to study infinitesimal
deformations and extensions of $Q$-bundles and this has been done on a
case by case basis.  For Bradlow pairs, see Thaddeus
\cite{thaddeus:1994} and for Higgs bundles see Nitsure
\cite{nitsure:1991}, Markman \cite{markman:1994}, Bottacin
\cite{bottacin:1995} and Hausel and Thaddeus
\cite{hausel-thaddeus:2000a}.  The study of extensions of triples is
important for the study of representations of the fundamental group of
a surface in the group $\mathrm{U}(p,q)$
\cite{bradlow-garcia-prada-gothen:2002,bradlow-garcia-prada-gothen:2003}.

The goal of this paper is to give a unified treatment of the homological
algebra of twisted quiver bundles, in order to study their
extensions and deformations.
The natural setting is to work in an abelian category with enough injectives,
but the category of $Q$-bundles is not abelian, because the category
of vector bundles over $X$ is itself not abelian.
We therefore work with representations of $Q$ 
in the category of sheaves of $\OO_X$-modules on a ringed space $(X,\OO_X)$,
which is an abelian category and has enough injectives. 
In this general setting, our results apply, for example,
both when $X$ is an algebraic variety with the
Zariski topology and $\OO_{X}$ is the sheaf of algebraic functions,
and also when $X$ is a complex manifold with the Hausdorff topology and
$\OO_{X}$ is the sheaf of holomorphic functions. We shall call a
representation of a quiver $Q$ in this category a \emph{$Q$-sheaf of
  $\OO_X$-modules} (or simply a \emph{$Q$-sheaf}, for short).  One
sees easily that the category of $Q$-sheaves is abelian.  In analogy
to the case of $Q$-vector spaces, the category of $Q$-sheaves of
$\OO_X$-modules is equivalent to the category of sheaves of modules
for a certain sheaf of algebras on $X$ (see
Section~\ref{sec:Q-sheaves}).  This more algebraic point of view is
important to our approach.

As in the vector bundle situation, the notion of $Q$-sheaf generalizes
to the notion of \emph{twisted $Q$-sheaf}; the twisting data consists
of a locally free sheaf $M_a$ of $\OO_X$-modules for each arrow $a \in
Q_1$ (again we refer to Section~\ref{sec:Q-sheaves} for precise
definitions).

Our main results are as follows.  We construct an explicit injective
resolution of a twisted $Q$-sheaf
(Proposition~\ref{prop:injective-resolution}); in particular it
follows that the abelian category of twisted $Q$-sheaves has enough
injectives.  This is done by combining the standard injective
resolution in the category of twisted $Q$-vector spaces with injective
resolutions in the category of sheaves.  From the explicit form of the
injective resolution of a twisted $Q$-sheaf we then deduce a long
exact sequence (Theorem~\ref{thm:main}) which relates the $\Ext$
groups in the category of twisted $Q$-sheaves to the $\Ext$ groups of
the sheaves which make up the twisted $Q$-sheaves.  This should be a
very useful tool in practical calculations and generalizes earlier
results for special quivers; it was probably first observed by Nitsure
\cite{nitsure:1991}, in the case of Higgs bundles.  Next we prove that
if $V$ and $W$ are twisted $Q$-sheaves such that $V$ is locally free, then the
group $\Ext^{i}(V,W)$ is isomorphic to a certain $i$th hypercohomology
group (Theorem~\ref{thm:ext=HH}).  This should be useful for
representing extension classes using $\mathrm{\check{C}}$ech or
Dolbeault cohomology, and also connects our results to the special
cases mentioned above, which are generally stated using
hypercohomology.

To finish this introduction, we mention some related work.  On the
subject of twisted $Q$-bundles, \'Alvarez-C\'onsul and
Garc{\'\i}a-Prada
\cite{alvarez-consul-garcia-prada:2001a,alvarez-consul-garcia-prada:2001b}
have recently proved a Hitchin-Kobayashi correspondence between stable
twisted $Q$-bundles and solutions to certain gauge theory equations
for special metrics on the underlying smooth bundles.  (In fact they
consider the even more general notion of quivers with relations).
This generalizes earlier results in the special cases mentioned
earlier.

We also should mention the paper \cite{biswas-ramanan:1994} by Biswas
and Ramanan.  They study the infinitesimal deformation space of
principal Higgs bundles on a curve using hypercohomology.  Even though
this is a somewhat different situation from ours, there are certain
similarities: for example they show that there is a long exact
sequence relating the infinitesimal deformation space to certain sheaf
cohomology groups of vector bundles associated to the underlying
principal bundle.

\paragraph{Acknowledgements.}
The first author thanks The Institute for Mathematical Sciences, 
University of Aarhus for hospitality while part of this work 
was carried out.  

\section{The path algebra} \label{sec:path-algebra}

In this section we recall a few relevant facts about representations
of quivers and the corresponding notion of modules over the path
algebra.  For more details see, e.g.,
\cite{auslander-reiten-smalo:1995}, \cite{benson:1998},
\cite{crawley-boevey:1992} or \cite{gabriel-roiter:1992},
but note that we use a slight modification of the standard definition
of path algebra.

Let $Q$ be a quiver.  A \emph{non-trivial path} in $Q$ is a sequence
$p = a_m \cdots a_0$ of arrows $a_i$, such that $ha_{i-1} = ta_{i}$,
for each $i$.
  \begin{displaymath}
    \xymatrix{
      {\bullet} & {\bullet} \ar[l]^{a_m} & \cdots \ar[l]^{a_{m-1}}
      & {\bullet} \ar[l]^{a_1}
      & {\bullet} \ar[l]^{a_0}
      }
  \end{displaymath}
We write $tp = ta_{0}$ and $hp = ha_{m}$.
In addition, for each vertex $i \in Q_0$, there is a
\emph{trivial path} $\langle i \rangle$, which 
begins and ends at $i$.
The \emph{length} of a non-trivial path is the number of arrows 
in the sequence; a trivial path has length zero.

Now, for a commutative ring $R$,
suppose that we have an $R$-module $M_{a}$ for each 
$a \in Q_{1}$.  
For each non-trivial path $p = a_m \cdots a_0$, we set
\begin{displaymath}
 M_{p} = M_{a_{m}} \otimes_{R} \cdots \otimes_{R} M_{a_{0}}.
\end{displaymath}
For each trivial path $\langle i \rangle$, we set 
$M_{\langle i \rangle} = R$ and write 
$e_i\in M_{\langle i \rangle}$ for the
element thereby identified with $1\in R$.

\begin{definition}
  \label{def:twisted-path-algebra}
    Let $Q$ be a quiver and let $M = \{M_{a}\}_{a \in Q_{1}}$ be a
    collection of $R$-modules.
    The \emph{$M$-twisted path algebra over $R$} is
    \begin{displaymath}
      A = \bigoplus_{p} M_{p}\ ,
    \end{displaymath}
    where the sum is over all paths in $Q$ and the $R$-linear
    multiplication is defined by
    \begin{displaymath}
      x_{p}x_{q} = x_{p} \otimes x_{q} 
      \ \text{if $tp=hq$, and $0$ otherwise,} 
    \end{displaymath}
    for $x_{p} \in M_{p}$ and $x_{q} \in M_{q}$, while
    \begin{align*}
      x_{p}e_{i} &= x_{p} 
      \ \text{if $i=tp$, and $0$ otherwise,} \\
      e_{i}x_{p} &= x_{p}
      \ \text{if $i=hp$, and $0$ otherwise,} \\
      e_ie_j &= e_i
      \ \text{if $i=j$, and $0$ otherwise.}
    \end{align*}    
\end{definition}

\begin{remark}
    (1) The $R$-algebra $A$ is graded, 
       by the subspaces corresponding to paths of length $n$, that is 
\[ 
  A_{n} = \bigoplus_{|p|=n} M_{p}\ .
\]
  Indeed, $A$ may also be defined as the tensor algebra, 
  over the ring $A_0=\bigoplus_{i\in Q_0} R$,
  of the $A_0$-bimodule $A_1= \bigoplus_{a\in Q_1} M_{a}$. 

    (2) When each $M_{a}=R$ (or is free with one generator),
    then we get the usual path algebra of the quiver.
    
    (3) Suppose that $a$ and $b$ are arrows with the same head and
    tail.  If we replace them by one arrow $c$ and let $M_{c} = M_{a}
    \oplus M_{b}$, then we clearly get the same algebra.
\end{remark}

The notion of a representation of $Q$ in the category of $R$-modules
defined in the Introduction can be generalized to take into account
twisting by the modules $M_{a}$ as follows: an \emph{$M$-twisted
  representation of $Q$} consists of $R$-modules $V_{i}$ for $i \in
Q_{0}$ and $R$-module maps
$$
\phi_{a} \colon M_a \otimes_{R} V_{ta} \lto V_{ha}
$$
for $a \in Q_1$.  A
\emph{morphism} between twisted representations $(V,\phi)$ and 
$(W,\psi)$ is a
collection of $R$-module maps $f_i \colon V_i \lto W_i$ such that the
diagrams
\begin{displaymath}
  \xymatrix{
    M_a \otimes_{R} V_{ta}
    \ar[d]^{1 \otimes f_{ta}} \ar[r]^-{\phi_{a}} & V_{ha} \ar[d]^{f_{ha}} \\
    M_a \otimes_{R} W_{ta} \ar[r]^-{\psi_{a}}                 & W_{ha} \\
  }
\end{displaymath}
commute for all $a\in Q_1$.
It follows immediately that the forgetful functors
$(V,\phi)\mapsto V_i$ are exact, for each $i\in Q_0$.

It is a straightforward exercise to show directly that the 
category of $M$-twisted
representations of $Q$ is an abelian category, using the fact that
the category of $R$-modules is abelian and the functor $M_a\otimes_{R}-$
is right exact.
However, this is also evident from the following.

\begin{proposition}
  \label{prop:category-equivalence}
  The category of twisted representations of $Q$ is equivalent to the
  category of left $A$-modules.  \qed
\end{proposition}

\begin{proof}
Let $(V,\phi)$ be a twisted representation of $Q$.  We can then define a left
$A$-module $V$ as follows: as an $R$-module it is
\begin{displaymath}
  V = \bigoplus_{i \in Q_0} V_i,
\end{displaymath}
and the $A$-module structure is extended linearly from 
\begin{align*}
  e_i v &= 
  \begin{cases}
    v, & v \in V_i, \\
    0, & v \in V_j \quad\text{for $j\neq i$,}
  \end{cases} \\
\intertext{for $i \in Q_{0}$ and}
  x_a v &= 
  \begin{cases}
    \phi_{a}(x_{a} \otimes v), & v \in V_{ta}, \\
    0, &v \in V_j \quad\text{for $j\neq ta$,}
  \end{cases}
\end{align*}
for $x_{a} \in M_{a}$ and $a \in Q_{1}$.  This construction can be
inverted as follows: given a left $A$-module $V$ we set
\begin{math}
  V_i = e_i V
\end{math}
for $i \in Q_0$ and define the map $\phi_{a} \colon M_a \otimes_{R} V_{ta}
\lto V_{ha}$ by
\begin{math}
  x_a \otimes v \lmapsto x_a v
\end{math}.
One quickly checks that morphisms of representations of $Q$ correspond to
$A$-module homomorphisms (see \cite{auslander-reiten-smalo:1995} or
\cite{benson:1998} for details).
\end{proof}

We will exploit this result and henceforth describe a twisted representation
$(V,\phi)$ of $Q$ as just an $A$-module $V$. 
Note in particular that we have a canonical identification $V_i = e_i V$.

The following $A$-module resolution is of fundamental importance to us. 

\begin{proposition}
\label{prop:twisted-exact-sequence}
Let $V$ be a module over the twisted path algebra $A$.  There
is an exact sequence
\begin{multline}
  \label{eq:twisted-exact-sequence}
  0 \lto V
  \overset{\epsilon}{\lto}
  \bigoplus_{i \in Q_0} \Hom_{R}(e_i A, V_i) \\
  \overset{d}{\lto}
  \bigoplus_{a \in Q_1}
  \Hom_{R}(M_a \otimes_{R} e_{ta}A, V_{ha})
  \lto 0, 
\end{multline}
where
\begin{displaymath}
  \epsilon(v)(x) = x v
\end{displaymath}
and the map $d$ is given by
\begin{displaymath}
  d \colon
  (\alpha_i)_{i\in Q_0} \mapsto
  \bigl( \alpha_{ha}\circ \mu_a - 
    \phi_a\circ (1\otimes \alpha_{ta}) \bigr)_{a\in Q_1}
\end{displaymath}
where $\mu_a : M_a\otimes e_{ta}A \to e_{ha}A$ is the multiplication map.
\end{proposition}

\begin{proof}
  Note that there is a natural inclusion 
  \begin{displaymath}
    \Hom_{A}(A,V) \into \bigoplus_{i\in Q_0} \Hom_{R}(e_i A, V_i)
  \end{displaymath}
  and that, under this inclusion, $\epsilon$ is just the standard
  isomorphism $V \xrightarrow{\cong}\Hom_{A}(A,V)$.  In particular,
  $\epsilon$ is injective.

  Next we show that $\im(\epsilon) = \ker(d)$.  Let
  $$
  \alpha = (\alpha_i)_{i\in Q_0}
  \in \bigoplus_{i\in Q_0} \Hom_{R}(e_i A, V_i)
  $$
  be such that $d\alpha = 0$.  The map $d$ can be written as 
  $d = (d_a)_{a \in Q_1}$ with 
  \begin{displaymath}
    d_a(\alpha)(x_a \otimes x) = \alpha_{ha}(x_a x) - x_a 
    \alpha_{ta}(x),
  \end{displaymath}
  and thus $d\alpha = 0$ means that for any $x \in V_{ta}$  and
  $x_a \in M_a$ we have
  \begin{displaymath}
    \alpha_{ha}(x_ax) = x_a \alpha_{ta} (x).
  \end{displaymath}
  Hence, we see that $d\alpha = 0$ if and only if $\alpha(yx) =
  y\alpha(x)$ for any $y \in A$; in other words if and only if
  $\alpha$ is $A$-linear.  Therefore $\ker(d) = \Hom_{A}(A,V)$, showing
  that $\im(\epsilon) = \ker(d)$.

  Finally, we must show that $d$ is surjective.  Let
  \begin{displaymath}
    \beta = (\beta_a)_{a \in Q_1}
    \in \bigoplus_{a \in Q_1}
    \Hom_{R}(M_a \otimes_{R} e_{ta}A, V_{ha}),
  \end{displaymath}
  then it suffices to find $\alpha$ such that 
  \begin{math}
    d_a(\alpha) = \beta_a
  \end{math}
  for each $a \in Q_1$.  Note that the grading of $A$ gives
  a natural grading 
  \begin{math}
    e_{i} A = \bigoplus_{l \geq 0} e_{i} A_l
  \end{math},
  and that we have
  \begin{displaymath}
    e_{i}A_{l} = \bigoplus_{\{a \;:\; ha = i\}}
      M_a \otimes_{R} e_{ta}A_{l-1}.
  \end{displaymath}
  We define $\alpha = (\alpha_i)_{i\in Q_0}$ by induction on $l$ as follows:
  \begin{itemize}
    \item On $e_{i} A_0$, we set $\alpha_i = 0$ for all $i \in Q_0$. 
    \item Assume that $\alpha$ is defined on $e_{i} A_{l-1}$ for all
    $i$, then we define $\alpha_i$ on $e_{i} A_l$ by setting
    \begin{displaymath}
      \alpha_{i}(x_a \otimes x)
      = x_a \alpha_{ta}(x) + \beta_a(x_a \otimes x)
    \end{displaymath}
    for $x_a \otimes x \in M_a \otimes_{R} e_{ta} A_{l-1}$.
  \end{itemize}
  It is then easy to verify that $d_a(\alpha) =\beta_a$.
\end{proof}

Let $A$ be any $R$-algebra,
let $V$ be a left $A$-module, let $W$ be a right $A$-module and let
$L$ be an $R$-module.  There is the standard adjunction
\begin{equation}
    \label{eq:standard-adjunction}
    \Hom_{A}(V,\Hom_{R}(W, L)) \cong \Hom_{R}(W \otimes_{A} V, L).
\end{equation}
The isomorphism relates $f$ on the left and $g$ on the right 
when
\begin{displaymath}
  f(v)(w)= g(w \otimes v)
\end{displaymath}
for all $v \in V$ and $w \in W$.
This adjunction has several important consequences, 
and we will use similar results later in the paper.
Firstly, since
\begin{displaymath}
  V_i \cong e_i A \otimes_{A} V,
\end{displaymath}
it follows that for any left $A$-module $V$ and for any $R$-modules
$N$ and $L$ we have the adjunction formula
\begin{equation}
  \label{eq:adjunction-A}
  \Hom_A(V,\Hom_{R}(N \otimes_{R} e_iA,L))
  \cong \Hom_{R}(N \otimes_{R} V_i,L),
\end{equation}
and from this we deduce the following Proposition.
\begin{proposition}
  \label{prop:injective-A}
  Let $N$ be a free $R$ module and $L$ be an injective $R$ module.
Then the left $A$-module
    \begin{displaymath}
      \Hom_{R}(N \otimes_{R} e_i A, L)
    \end{displaymath}
is injective for each $i \in Q_0$.
\end{proposition}

\begin{proof}
  We need to see that the functor
\begin{displaymath}
  \Hom_A(\blank,\Hom_{R}(N \otimes_{R} e_iA,L))
\end{displaymath}
is exact, which it is, because, by \eqref{eq:adjunction-A}, it is the
composite of three exact functors: 
the forgetful functor $V \mapsto V_i$ is always exact,
$N \otimes_{R} \blank$ is exact when $N$ is free, and
$\Hom_{R}(\blank,L)$ is exact when $L$ is injective.
\end{proof}

As a special case, suppose that $R = \kk$ is a field. 
Then every vector space is both free and injective
and the sequence
\eqref{eq:twisted-exact-sequence} is the `standard'
injective resolution of $V$.

\begin{example}
  Consider the quiver \eqref{eq:higgs-quiver}.  In this case the path
  algebra $A$ can be identified with the polynomial algebra $\kk[x]$.
  Let $V$ be a module over $\kk [x]$, then the injective resolution
  \eqref{eq:twisted-exact-sequence} takes the form
  \begin{displaymath}
    0 \lto V \overset{\epsilon}{\lto} \Hom_{\kk}(\kk[x], V)
    \overset{d}{\lto}
    \Hom_{\kk}(\kk[x], V) \lto 0,
  \end{displaymath}
  where
  $\epsilon(v)(p)=pv$ and 
  $d(\alpha)(p) = \alpha(xp) - x \alpha(p)$.
\end{example}

\section{Quiver sheaves} \label{sec:Q-sheaves}

Let $Q$ be a quiver and let $(X,\OO_X)$ be a ringed space. In this
section, we consider twisted representations of $Q$ in the category
of sheaves of $\OO_X$-modules on $X$.
For basic definitions about ringed spaces and sheaves of
$\OO_X$-modules, refer to Hartshorne \cite{hartshorne:1977},
Ch.~II, Sec.~1,2,5.
Note also that the category of sheaves of $\OO_X$-modules on $X$ 
has enough injectives (\cite{hartshorne:1977} Prop.~III.2.2).

Suppose we are given a collection $\{M_a\}_{a\in Q_{1}}$ of locally free
sheaves of $\OO_X$-modules.  
An \emph{$M$-twisted $Q$-sheaf of $\OO_X$-modules}
(or a \emph{twisted $Q$-sheaf on $X$}, for short) consists of
$\OO_X$-modules $V_i$ for $i \in Q_0$ and morphisms
$$
\phi_{a} \colon M_a \otimes_{\OO_X} V_{ta} \lto V_{ha}
$$
for $a \in Q_1$.  A \emph{morphism} between twisted  $Q$-sheaves 
$(V,\phi)$ and $(W,\psi)$
on $X$ is a collection of $\OO_X$-module maps $f_i \colon V_i \lto
W_i$ such that the diagrams
\begin{displaymath}
  \xymatrix{
    M_a \otimes_{\OO_X} V_{ta}
    \ar[d]^{1 \otimes f_{ta}} \ar[r]^-{\phi_{a}} & V_{ha} \ar[d]^{f_{ha}} \\
    M_a \otimes_{\OO_X} W_{ta} \ar[r]^-{\phi_{a}}                 & W_{ha} \\
  }
\end{displaymath}
commute for all $a\in Q_1$.
As in Section~\ref{sec:path-algebra}, the forgetful functors
$(V,\phi)\mapsto V_i$ are exact, for each $i\in Q_0$,
and the category of $M$-twisted $Q$-sheaf of $\OO_X$-modules
is abelian.

Next we describe the corresponding algebraic notion,
analogous to the definition of the twisted path algebra in 
Section~\ref{sec:path-algebra}.  For each non-trivial path 
$p = a_{m} \cdots a_{0}$ in $Q$ we set 
$$
M_{p} = M_{a_{m}} \otimes_{\OO_{X}} \cdots \otimes_{\OO_{X}} M_{a_{0}}
$$
and we set $M_{\langle i \rangle} = \OO_{X}$.  We define the
\emph{$M$-twisted path algebra on $X$} to be the
sheaf of $\OO_X$-algebras $B = \bigoplus_{p} M_{p}$ where the sum is 
over all paths in $Q$ and the multiplication is the obvious one (cf.\
Definition~\ref{def:twisted-path-algebra}).

In a similar way to Proposition~\ref{prop:category-equivalence},
the category of twisted $Q$-sheaves is equivalent 
to the category of sheaves of $B$-modules on $X$ 
(a detailed proof can be found in
\cite{alvarez-consul-garcia-prada:2001b}).  
Thus we shall interchangeably use the expressions 
``twisted  $Q$-sheaf'' and ``sheaf of $B$-modules'', 
or, for short, ``$B$-module''.

Note that, for the elementary facts above,
we do not use the assumption that the twisting sheaves $M_a$ are locally free,
i.e. that the functor $M_a \otimes_{\OO_{X}}\blank$ is exact.
However, this assumption is used in realising 
the main objective of this section, that is, to show that
the category of twisted $Q$-sheaves has enough injectives.

Specifically, we shall 
construct an injective resolution of a
twisted $Q$-sheaf by combining the resolution
\eqref{eq:twisted-exact-sequence} in the $A$-module direction with an
injective resolution in the $\OO_{X}$-module direction.  
The first step is to generalize 
Proposition~\ref{prop:twisted-exact-sequence} to $B$-modules.

\begin{proposition}
\label{prop:sheaf-resolution}
Let $V$ be a $B$-module.  There is an exact sequence of sheaves
\begin{multline}
  \label{eq:twisted-sheaf-exact-complex}
  0 \lto V
  \overset{\epsilon}{\lto}
  \bigoplus_{i \in Q_0} \shom_{\OO_{X}}(e_i B, V_i) \\
  \overset{d}{\lto}
  \bigoplus_{a \in Q_1}
  \shom_{\OO_{X}}(M_a \otimes_{\OO_{X}} e_{ta}B, V_{ha})
  \lto 0, 
\end{multline}
where the map
$\epsilon$ is defined by
\begin{equation}
  \epsilon(v)(x) = x v \label{eq:sheaf-epsilon}
\end{equation}
for local sections $v$ and $x$ and the map $d$ is defined on a local 
section $\alpha = (\alpha_i)_{i\in Q_0}$ by
\begin{equation}
  d \colon
  (\alpha_i)_{i\in Q_0} \mapsto
  \bigl( \alpha_{ha}\circ \mu_a - 
    \phi_a\circ (1\otimes \alpha_{ta}) \bigr)_{a\in Q_1} 
    \label{eq:sheaf-d}
\end{equation}
where $\mu_a : M_a\otimes e_{ta}B \to e_{ha}B$ is the
multiplication map.
\end{proposition}

\begin{proof}
  It is sufficient to show exactness on stalks.  Let $\OO_x$ be the
  local ring at $x$, and let $M_{a,x}$, $B_x$ and $V_x$ be the stalks
  at $x$ of $M_a$, $B$ and $V$, respectively.  Recall from Serre
  \cite[n$^{\mathrm{o}}$ 10]{serre:1955} that the stalk of the tensor
  product of sheaves is the tensor product of the stalks over the
  local ring.  Thus we need to show that the sequence
\begin{multline*}
  0 \lto V_x
  \overset{\epsilon}{\lto}
  \bigoplus_{i \in Q_0} \Hom_{\OO_{x}}(e_i B_x, V_{i,x}) \\
  \overset{d}{\lto}
  \bigoplus_{a \in Q_1}
  \Hom_{\OO_{x}}(M_{a,x} \otimes_{\OO_{x}} e_{ta}B_x, V_{x,ha})
  \lto 0, 
\end{multline*}
is exact for every $x \in X$.
But, again using the fact on tensor products mentioned above, we see that
$B_x$ is the $M_{a,x}$-twisted path algebra over $\OO_x$.  Hence
exactness is immediate from
Proposition~\ref{prop:twisted-exact-sequence}.
\end{proof}

Another important ingredient in our construction is the following
analogue for sheaves of the adjunction formula
\eqref{eq:standard-adjunction}.

\begin{proposition}
  \label{prop:sheaf-adjunction}
  Let $B$ be a sheaf of $\OO_X$-algebras, let $V$ a sheaf of left
  $B$-modules, let $W$ be a sheaf of right $B$-modules and let $L$ be
  an $\OO_X$-module.  Then there is a natural isomorphism
  \begin{displaymath}
    \shom_{B}(V,\shom_{\OO_X}(W, L))
    \cong \shom_{\OO_X}(W \otimes_{B} V, L).    
  \end{displaymath}
\end{proposition}
\begin{proof}
  This follows from \eqref{eq:standard-adjunction} by an argument 
  analogous to the one used in the proof of 
  Proposition~\ref{prop:sheaf-resolution} above.
\end{proof}

\begin{corollary}
  \label{prop:adjunction}
  Let $V$ be a $B$-module and $N$ and $L$ be a $\OO_{X}$-modules. 
  Then there is a natural isomorphism
  \begin{displaymath}
    \shom_{\OO_{X}}(N \otimes_{\OO_X} V_i, L)
    \cong \shom_B(V,\shom_{\OO_{X}}(N \otimes_{\OO_X} e_i B, L)).
  \end{displaymath}
\end{corollary}

\begin{proof}
  This follows from Proposition~\ref{prop:sheaf-adjunction},
  because $V_i \cong e_i B \otimes_B V$.
\end{proof}

The following is analogous to Proposition~\ref{prop:injective-A}. 

\begin{proposition}
  \label{prop:injective-B-modules}
  Let $L$ be an injective $\OO_{X}$-module and let $N$ be a locally
  free $\OO_X$-module.  Then
\begin{displaymath}
  \shom_{\OO_{X}} (N \otimes_{\OO_X} e_{i} B, L)
\end{displaymath}
  is an injective $B$-module.
\end{proposition}

\begin{proof}
  We have to show that
\begin{displaymath}
  \Hom_B(\blank,\shom_{\OO_{X}} (N \otimes_{\OO_X} e_{i} B, L))
\end{displaymath}
is an exact functor.  By Corollary~\ref{prop:adjunction}
this functor is the composite of the three functors, $V \lmapsto
V_i$, $N \otimes_{\OO_{X}} \blank$ and $\Hom_{\OO_{X}}(\blank,L)$.
Just as in the proof of Proposition~\ref{prop:injective-A}, 
the first one is always exact 
the second one is exact because $N$ is locally free
and the third one is exact because $L$ is injective.
\end{proof}

Now, if $V$ is any $B$-module, then for each $\OO_{X}$-module $V_i$, 
we can choose an injective resolution 
\begin{displaymath}
  V_i \overset{\epsilon}{\lto} V_i^0
  \overset{\delta}{\lto} V_i^1 \overset{\delta}{\lto} \cdots
\end{displaymath}
Because $M_a \otimes_{\OO_{X}}\blank$ is an exact functor,
the usual argument, using injectivity of the $V_{ha}^n$,
shows that the maps
\begin{equation*}
  \phi_{a} \colon M_a \otimes_{\OO_{X}} V_{ta} 
  \lto V_{ha}
\end{equation*}
lift to chain maps
\begin{equation}
 \label{eq:chain-lift}
  \phi^{\bullet}_{a} \colon M_a \otimes_{\OO_{X}} V_{ta}^{\bullet} 
  \lto V_{ha}^{\bullet}
\end{equation}
for each $a \in Q_1$.
Thus the $V^{n} = \bigoplus_{i} V_i^{n}$ become
$B$-modules and the maps $\delta \colon V^{n} \to V^{n+1}$ become
$B$-module homomorphisms.

We construct the desired injective resolution of
$V$ from the following double complex, which we denote by
$C(V)^{\bullet\bullet}$.
\begin{equation}
  \label{eq:twisted-double-complex}
  \xymatrix{
    {\vdots} & {\vdots} \\
    {\bigoplus_{i \in Q_0} \shom_{\OO_{X}}(e_i B, V_i^1)}
        \ar[r]\ar[u]
        & {\bigoplus_{a \in Q_1}
          \shom_{\OO_{X}}(M_a \otimes_{\OO_{X}} e_{ta} B, V_{ha}^1)}
        \ar[u] \\
    {\bigoplus_{i \in Q_0} \shom_{\OO_{X}}(e_i B, V_i^0)}
        \ar[r]\ar[u]
        & {\bigoplus_{a \in Q_1}
          \shom_{\OO_{X}}(M_a \otimes_{\OO_{X}} e_{ta} B, V_{ha}^0)}
        \ar[u] \\
  }
\end{equation}
The horizontal differentials are $(-1)^{q} d$, where $d$ is the 
map \eqref{eq:sheaf-d} for $V^{q}$,
and the vertical differentials are induced by
$\delta \colon V_i^{n} \lto V_i^{n+1}$.
Using the fact that the $\delta$ are $B$-module homomorphisms it is
easy to check that the squares of the double complex anti-commute.

Denote the total complex of $C(V)^{\bullet\bullet}$ by
$I(V)^{\bullet}$ and note that we have a canonical injection of
$B$-modules
\begin{displaymath}
  \bar{\epsilon} \colon V \lto I(V)^{0}
    = \bigoplus_{i \in Q_0} \shom_{\OO_{X}}(e_i B, V_i^0)
\end{displaymath}
where $\bar{\epsilon}$ is 
given by composing 
\begin{displaymath}
  \epsilon \colon V^{0}
  \lto
  \bigoplus_{i \in Q_0} \shom_{\OO_{X}}(e_i B, V_i^{0})
\end{displaymath}
defined in \eqref{eq:sheaf-epsilon} with the inclusion $V \into V^{0}$.

\begin{proposition}
  \label{prop:injective-resolution}
  For any $B$-module $V$, 
  the complex $I(V)^{\bullet}$ is an injective resolution of $V$.
\end{proposition}

\begin{proof}
  It follows from Proposition~\ref{prop:injective-B-modules} that
  each term in \eqref{eq:twisted-double-complex} is injective
  and so $I(V)^{n}$ is an injective $B$-module for each $n$. 
  We have just seen that $V$ injects into $I(V)^{0}$,
  so all that remains to check is that 
  $V \to I(V)^{\bullet}$ is an exact complex.
  We use the following notation:
  \begin{align*}
    F(V) &= \bigoplus_{i \in Q_0} \shom_{\OO_{X}}(e_i B, V_i)\ , \\
    G(V) &= \bigoplus_{a \in Q_1}
     \shom_{\OO_{X}}(M_a \otimes_{\OO_{X}} e_{ta}B, V_{ha})\ .
  \end{align*}
  Then we have the following commutative diagram:
  \begin{displaymath}
  \xymatrix{
   & & {\vdots} & {\vdots} & \\
   & & C(V)^{01} \ar[r]\ar[u] & C(V)^{11} \ar[u] & \\
   & & C(V)^{00} \ar[r]\ar[u] & C(V)^{10} \ar[u] & \\
   0 \ar[r] & V \ar[ur]^{\bar{\epsilon}} \ar[r]^{\epsilon}
     & F(V) \ar[r]^{d}\ar[u] & G(V) \ar[r]\ar[u] & 0
  }
  \end{displaymath}
  Since, for any $i \in Q_0$ and $a \in Q_1$, the $\OO_X$-modules
  $e_iB$ and $M_a \otimes_{\OO_X} e_{ta} B$ are locally free, we see
  that $\shom_{\OO_X}(e_i B, \blank)$ and
  $\shom_{\OO_X}(M_a\otimes_{\OO_X} e_i B, \blank)$ are exact
  functors.  Hence the columns of the above diagram are exact.
  Furthermore, the bottom row of the diagram is exact from
  Proposition~\ref{prop:sheaf-resolution}.  From these two facts it
  is easy to see that $V \to I(V)^{\bullet}$ is an exact complex.
\end{proof}

\section{A long exact sequence} \label{sec:main-theorem}

Now that we can take injective resolutions of $B$-modules,
we may define $\Ext_B^i(V,W)$
as the $i$th cohomology of the complex
$\Hom_{B}(V,I(W)^{\bullet})$,
where $I(W)^{\bullet}$ is an injective resolution of $W$.  
Our main theorem provides a tool for calculating
$\Ext^i_B(V,W)$, in terms of the $\Ext$ groups of the 
$\OO_X$-modules $V_{i}$ and $W_{j}$, via a long exact sequence.
As explained in the Introduction, several earlier results are special
cases of this.

Before stating the result, we need to explain how to interpret 
some of the notation.
We are already familiar with the map
\[
 \Hom(E,F) \to \Hom(M_a\otimes E,M_a\otimes F) : f\mapsto 1\otimes f
\]
given by the functor $E\mapsto M_a\otimes E$.
Since $M_a$ is locally free, 
this functor is exact and hence it also induces a map
\[
 \Ext^i(E,F) \to \Ext^i(M_a\otimes E,M_a\otimes F),
\]
which we shall also denote by $f\mapsto 1\otimes f$.
This map is easily understood in the Yoneda picture,
where $\Ext^i$ 
classifies exact sequences
from $F$ to $E$ with $i$ intermediate terms.
However, for our purposes we need the following equivalent 
formulation.
If $F^{\bullet}$ is any resolution of $F$
and $\ftil\in \Hom(E,F^i)$ is a cocyle representing 
$f\in \Ext^i(E,F)$, then 
$1\otimes \ftil\in \Hom(M_a\otimes E,M_a\otimes F^i)$
is a cocyle representing $1\otimes f\in\Ext^i(M_a\otimes E,M_a\otimes F)$.

With this understood, we can state and prove our theorem.
The context is the general one of Section~\ref{sec:Q-sheaves}:
$(X,\OO_X)$ is a ringed space, $Q$ is a quiver and $M_a$ are
locally free $\OO_X$-module for each $a \in Q_1$.
Then $B$ is the $M$-twisted path algebra over $X$.

\begin{theorem}
  \label{thm:main}
  For any $B$-modules $V$ and $W$, there is a long exact sequence
  \begin{eqnarray*}
    0 \to \Hom_B(V,W)
    &\to& \bigoplus_{i\in Q_0} \Hom_{\OO_X}(V_{i},W_{i}) 
   \overset{\partial}{\lto}   
   \bigoplus_{a\in Q_1} \Hom_{\OO_X}
      (M_{a} \otimes V_{ta},W_{ha}) \\
    \to \Ext^1_B(V,W) 
    &\to& \bigoplus_{i\in Q_0} \Ext^1_{\OO_X}(V_{i},W_{i})
   \overset{\partial}{\lto}   
   \bigoplus_{a\in Q_1} \Ext^1_{\OO_X}
      (M_{a} \otimes V_{ta},W_{ha}) \\
    \to \Ext^2_B(V,W) &\to& \cdots .
  \end{eqnarray*}
The maps $\partial$ are given by
\[
\partial\colon
 \bigl(f_i\bigr)_{i\in Q_0}  \lmapsto  
 \bigl( f_{ha}\circ \phi_a - \psi_a\circ (1\otimes f_{ta}) \bigr)_{a\in Q_1}
\]
where $\circ$ is the Yoneda product.

\end{theorem}

\begin{proof}
Let $I(W)^{\bullet}$ be the injective resolution of $W$ given by
Proposition~\ref{prop:injective-resolution},
i.e. the total complex of the double complex 
$C(W)^{\bullet\bullet}$ of \eqref{eq:twisted-double-complex}. 
Then the group $\Ext^i_B(V,W)$ is 
the $i$th cohomology group of the complex
$\Hom_B(V,I(W)^{\bullet})$,
which can be computed from the spectral sequence for
the double complex 
$\Hom_B(V,C(W)^{\bullet\bullet})$.
Using the adjunction of Corollary~\ref{prop:adjunction}, 
this double complex can be written
\begin{equation}
  \label{eq:OX-double-complex}
  \xymatrix{
    {\vdots} & {\vdots} \\
    {\bigoplus_{i \in Q_0} \Hom_{\OO_{X}}(V_i, W_i^1)}
        \ar[r]\ar[u]
        & {\bigoplus_{a \in Q_1}
          \Hom_{\OO_{X}}(M_a \otimes V_{ta}, W_{ha}^1)}
        \ar[u] \\
    {\bigoplus_{i \in Q_0} \Hom_{\OO_{X}}(V_i, W_i^0)}
        \ar[r]\ar[u]
        & {\bigoplus_{a \in Q_1}
          \Hom_{\OO_{X}}(M_a \otimes V_{ta}, W_{ha}^0)}
        \ar[u] \\
  }
\end{equation}
If we pass to the $E_1$ term of the spectral sequence by computing
vertical cohomology, then we obtain the required groups 
\[
 \textstyle{\bigoplus_{i \in Q_0}} 
 \Ext^q_{\OO_X}(V_{i},W_{i})
 \quad\text{and}\quad 
 \textstyle{\bigoplus_{a \in Q_1}}
 \Ext^q_{\OO_X}(M_{a} \otimes V_{ta},W_{ha}).
\]
The fact that we started with a two column complex means that
the spectral sequence converges at the $E_2$ term and yields
precisely the long exact sequence of the theorem.

To see that the maps are as claimed, observe that the horizontal maps
in \eqref{eq:OX-double-complex} are given by
\begin{equation}\label{eq:ftil}
\bigl(\ftil_i\bigr)_{i\in Q_0}  \lmapsto  
 \bigl( \ftil_{ha}\circ \phi_a 
- \psi^q_a\circ (1\otimes \ftil_{ta}) \bigr)_{a\in Q_1}
\end{equation}
where $\ftil_i\in\Hom(V_i,W_i^q)$ is a cocycle representing 
$f_i\in\Ext^q(V_i,W_i)$, 
while $\psi^q_a\in\Hom(M_{a} \otimes W_{ta}^q,W_{ha}^q)$ 
are the lifts of $\psi_a\in \Hom(M_{a} \otimes W_{ta},W_{ha})$, 
as in \eqref{eq:chain-lift}.
In \eqref{eq:ftil} the symbol $\circ$ just denotes composition
of maps.

The first composite 
$\ftil_{ha}\circ \phi_a$ is immediately seen to
be a cocyle representing the Yoneda product $f_{ha}\circ \phi_a$.
For the second composite, observe, 
from the remark preceding the statement of the theorem, 
that 
\[
1\otimes \ftil_{ta}\in \Hom(M_{a} \otimes V_{ta},M_{a} \otimes W_{ta}^q)
\]
is a cocyle representing 
$1\otimes f_{ta}\in\Ext^q(M_{a} \otimes V_{ta},M_{a} \otimes W_{ta})$
and thus that the composite $\psi^q_a \circ (1\otimes \ftil_{ta})$
is a cocyle representing the Yoneda product 
$\psi_a \circ (1\otimes f_{ta})$

Note for the careful reader:
because it makes no difference to the result,
for simplicity we have not carried the signs $(-1)^q$ attached to the
maps $d$ of \eqref{eq:twisted-double-complex}
through to the maps $\partial$ of the theorem.
\end{proof}

\section{Hypercohomology} \label{sec:hypercohomology}

Let $V$ and $W$ be $B$-modules.  When $V$ is locally free, we show in
this section how to give a description of $\Ext^{i}(V,W)$ in terms of
hypercohomology.

Consider the exact sequence of sheaves
given by Proposition~\ref{prop:sheaf-resolution}:
\begin{multline}
  \label{eq:sheaf-complex-W}
  0 \lto W
  \overset{\epsilon}{\lto}
  \bigoplus_{i \in Q_0} \shom_{\OO_{X}}(e_i B, W_i) \\
  \overset{d}{\lto}
  \bigoplus_{a \in Q_1}
  \shom_{\OO_{X}}(M_{a} \otimes_{\OO_{X}} e_{ta}B, W_{ha})
  \lto 0.
\end{multline}
Applying the functor $\shom_{B}(V,\blank)$ and using
Corollary~\ref{prop:adjunction} we obtain the exact sequence
\begin{multline}
  \label{eq:shom_B-sequence}
  0 \lto \shom_B(V,W)
  \lto
  \bigoplus_{i \in Q_0} \shom_{\OO_{X}}(V_i, W_i) \\
  \overset{\delta}{\lto}
  \bigoplus_{a \in Q_1}
  \shom_{\OO_{X}}(M_{a} \otimes_{\OO_{X}}V_{ta}, W_{ha}),
\end{multline}
where $\delta$ is given by
\begin{equation}
\label{eq:defD}
\bigl(f_i\bigr)_{i\in Q_0}  \lmapsto  
 \bigl( f_{ha}\circ \phi_a - \psi_a\circ (1\otimes f_{ta}) \bigr)_{a\in Q_1}
\end{equation}
The complex of sheaves 
\begin{displaymath}
\label{eq:sheaf-complex-2}
  C^{\bullet}(V,W) : \quad
  \bigoplus_{i \in Q_0} \shom_{\OO_X}(V_i, W_i)
  \overset{\delta}{\lto}
  \bigoplus_{a \in Q_1}
  \shom_{\OO_{X}}(M_{a} \otimes_{\OO_{X}}V_{ta}, W_{ha})
\end{displaymath}
formed by the last two terms of this sequence is important because of
the following result.
\begin{theorem}
  \label{thm:ext=HH}
  Let $(X,\OO_X)$ be a ringed space, let $Q$ be a quiver, let $M_a$ be
  a locally free $\OO_X$-module for each $a \in Q_1$ and let $B$ be
  the $M$-twisted path algebra over $X$.  Let $V$ and $W$ be
  $B$-modules and suppose that $V$ is locally free.  Then
  \begin{displaymath}
    \Ext_B^p(V,W) \cong \HH^p(C^{\bullet}(V,W)).
  \end{displaymath}
\end{theorem}

We shall need the following three lemmas for the proof.
\begin{lemma}\label{lem:H-gamma-C}
  Let $V$ and $W$ be $B$-modules and consider the two term complex
  obtained by applying the global sections functor $\Gamma$ to
  $C^{\bullet}(V,W)$.  Then
  \begin{align*}
    H^0(\Gamma(C^{\bullet}(V,W))) &\cong \Hom_B(V,W). \\
  \intertext{If, moreover, $W$ is injective, then}
    H^p(\Gamma(C^{\bullet}(V,W))) &= 0
  \end{align*}
  for $p \geq 1$.
\end{lemma}

\begin{proof}
  Observe that the complex 
  $\Gamma(C^{\bullet}(V,W))$ coincides with the part of the long exact 
  sequence of Theorem~\ref{thm:main} made up by the second and third 
  term.  The result follows.  
\end{proof}

\begin{lemma}\label{lem:W_i-injective}
  Let $W$ be an injective $B$-module.  Then each $W_{i}$ is an 
  injective $\OO_{X}$-module.
\end{lemma}

\begin{proof}
  We claim that the left adjoint of the forgetful functor $W \mapsto
  W_{i}$ from $B$-modules to $\OO_{X}$-modules is the functor $L
  \mapsto Be_{i} \otimes L$.  This latter functor is exact because the
  $M_{a}$ are locally free and hence its right
  adjoint (the forgetful functor) preserves injectives.  
  Thus the lemma follows from the claim.
  
  In order to prove the claim we use the adjunction
  \begin{displaymath}
    \Hom_{B}(V \otimes L,W) 
    \cong \Hom_{\OO_{X}}(L,\Hom_{B}(V,W))
  \end{displaymath}
  for an $\OO_{X}$-module $L$ and $B$-modules $V$ and $W$: taking $V = 
  Be_{i}$ it follows that the functor $L
  \mapsto Be_{i} \otimes L$ is the left adjoint of $W \mapsto 
  \Hom_{B}(Be_{i},W)$.  We finish the proof by observing that there is 
  a natural isomorphism $W_{i} = e_{i} W \overset{\cong}{\to} 
  \Hom_{B}(Be_{i},W)$ given by $e_{i}w \mapsto (xe_{i} \mapsto xw)$, 
  with inverse $\alpha \mapsto \alpha(e_{i})=e_{i}(\alpha(e_{i}))$.
\end{proof}

\begin{lemma}\label{lem:injective-OO}
  Let $L$ and $M$ be locally free $\OO_{X}$-modules and let $I$ be an 
  injective $\OO_{X}$-module.  Then $\shom_{\OO_{X}}(M \otimes L,I)$ 
  is an injective $\OO_{X}$-module.
\end{lemma}

\begin{proof}
  This follows from the standard adjunction 
  \begin{displaymath}
    \Hom_{\OO_{X}}(N,\shom_{\OO_{X}}(M \otimes L,I))
    \cong \Hom_{\OO_{X}}(M \otimes L \otimes N,I)
  \end{displaymath}
  by an argument similar to that used in the proof of 
  Propositions~\ref{prop:injective-A} and 
  \ref{prop:injective-B-modules}. 
\end{proof}

\begin{proof}[Proof of Theorem~\ref{thm:ext=HH}]
  We begin by noting that the statement is true for $p=0$.  For any
  complex of sheaves $\mathcal{F}^{\bullet}$ one has 
  \begin{math}
    \HH^0(\mathcal{F}^{\bullet}) = H^0(\Gamma(\mathcal{F}^{\bullet}))
  \end{math}
  and therefore we have, using Lemma~\ref{lem:H-gamma-C}, that 
  \begin{displaymath}
    \HH^0(C^{\bullet}(V,W)) \cong 
    H^0(\Gamma(C^{\bullet}(V,W)))\cong \Hom_B(V,W).
  \end{displaymath}
  
  Thus it suffices to show that the functor
  $\HH^p(C^{\bullet}(V,\blank))$ on $B$-modules is a universal
  $\delta$-functor.  Let $0 \to W' \to W \to W'' \to 0$ be a short
  exact sequence of $B$-modules.  Then, because $M_a$ and $V$ are
  locally free, we get a corresponding short exact sequence of
  complexes
  \begin{displaymath}
    0 \lto C^{\bullet}(V,W') \lto C^{\bullet}(V,W) \lto
    C^{\bullet}(V,W'') \lto 0.
  \end{displaymath}
  Hence, by general properties of hypercohomology, 
  we get the required long exact sequence 
  \begin{align*}
    0 &\lto \HH^0(C^{\bullet}(V,W')) \lto \HH^0(C^{\bullet}(V,W)) \lto
    \HH^0(C^{\bullet}(V,W'')) \\
    &\lto \HH^1(C^{\bullet}(V,W')) \lto
    \cdots\ ,
  \end{align*}
  natural with respect to morphisms of short exact sequences.  
  Thus $\HH^p(C^{\bullet}(V,\blank))$ is a $\delta$-functor. 
  
  It remains to show that $\HH^p(C^{\bullet}(V,\blank))$ is universal.
  Since the category of $B$-modules has enough injectives it suffices
  to show that for any injective $B$-module $I$, one has
  $\HH^p(C^{\bullet}(V,I)) = 0$ for all $p \geq 1$.
    Consider the first hypercohomology spectral sequence
  \begin{displaymath}
    {}^{I}E_{2}^{pq} = H^{p}(R^q\Gamma(C^{\bullet}(V,I))) 
    \implies
    \HH^{p+q}{C^{\bullet}(V,I)}.
  \end{displaymath}
    (Here we have denoted sheaf cohomology by $R^q\Gamma$ to avoid 
    confusion, i.e., $R^q\Gamma(C^{\bullet}(V,I))$ is the 
    complex obtained from $C^{\bullet}(V,I)$ by taking the $q$th 
    sheaf cohomology group of each term.)  By 
    Lemma~\ref{lem:W_i-injective} we have that each $I_{i}$ is an 
    injective $\OO_{X}$-module.  Hence
    Lemma~\ref{lem:injective-OO} shows that each term of 
    $C^{\bullet}(V,I)$ is an injective $\OO_{X}$-module and thus has 
    vanishing higher sheaf cohomology.  It follows that the spectral 
    sequence collapses at the $E_{2}$-term and we get
    \begin{displaymath}
      \HH^p(C^{\bullet}(V,I)) = 
      H^{p}(\Gamma(C^{\bullet}(V,I))) 
    \end{displaymath}
    which vanishes for $p \geq 1$ by Lemma~\ref{lem:H-gamma-C}.  This 
    finishes the proof.
\end{proof}


\providecommand{\bysame}{\leavevmode\hbox to3em{\hrulefill}\thinspace}

\vspace{1cm}

\newlength{\width}
\settowidth{\width}{\small Departamento de Matem{\'a}tica Pura}
\noindent
\parbox[t]{\width}{\small
  Departamento de Matem{\'a}tica Pura \\
  Faculdade de Ci{\^e}ncias \\
  Universidade do Porto \\
  Rua do Campo Alegre \\
  4169-007 Porto, Portugal \\
  E-mail: \texttt{pbgothen@fc.up.pt} \\
}
\hfill
\settowidth{\width}{\small E-mail: \texttt{A.D.King@maths.bath.ac.uk}}
\parbox[t]{\width}{\small
  Mathematical Sciences \\
  University of Bath \\
  Claverton Down \\
  Bath BA2 7AY, U.K. \\
  E-mail: \texttt{A.D.King@maths.bath.ac.uk} \\
}

\end{document}